\documentclass{amsart}
\usepackage{amssymb}
\usepackage{amsmath}
\usepackage{hyperref}

\usepackage{colortbl}
\usepackage{kotex}

\usepackage{tikz}

\newtheorem{thm}{Theorem}[section]
\newtheorem{lem}[thm]{Lemma}

\theoremstyle{definition}

\newtheorem{rem}[thm]{Remark}

\newtheorem{prop}[thm]{Proposition}

\theoremstyle{remark}

\numberwithin{equation}{section}

\begin{document}

\title[Four vertices of a square of area $n$]
{Finding four vertices of a square of area $n$ in a grid by geometric construction}

\author{Dohyeok Kim}
\address{Department of Mathematics Education,
Kangwon National University, Chuncheon 24341, Republic of Korea}
\email{ehgur135@gmail.com}

\author{Yunsoo Jang}

\address{Department of Mathematics Education,
Kangwon National University, Chuncheon 24341, Republic of Korea}
\email{yunsoojang@kangwon.ac.kr}


\thanks{Corresponding author: Yunsoo Jang, yunsoojang@kangwon.ac.kr}
\thanks{
D. Kim and Y. Jang were supported by Basic Science Research Program through the 
National Research Foundation of Korea(NRF) funded by the Ministry of Education(No. RS-2025-25415913).
}




\begin{abstract}
In this study, we consider geometric construction in an infinite square grid.
Our purpose is to find efficient methods for constructing a square of area $n\in\mathbb{N}$ based on arithmetic properties of the sum of squares.
As a consequence, for each $n\in\mathbb{N}$, we derive the optimal number of construction steps required to determine the four vertices of the square.
Furthermore, we classify all positive integers into three classes according to these optimal numbers.

\end{abstract}

\maketitle

\section{\bf Introduction}
\label{sec1}

We discuss the way to draw a square of area $n\in\mathbb{N}$ in a grid by geometric construction.
For this, we consider an infinite square grid in $\mathbb{R}^2$ with side length $1$. 
As a simple example, the following figure illustrates how to describe a square of area $2$ in this grid.

\begin{figure}[htbp]

\centering
\begin{tikzpicture}[scale=2]
    \draw[thin, gray!50] (0,0) grid (2,2);
    
    \fill[black, opacity=0.4] (1,0) -- (2,1) -- (1,2) -- (0,1) -- cycle;
    
    \filldraw[black] (1,0) circle (1.2pt)
                     (2,1) circle (1.2pt)
                     (1,2) circle (1.2pt)
                     (0,1) circle (1.2pt);
\end{tikzpicture}
\caption{An example of a square of area $2$ in a grid}
\label{fig1}
\end{figure}
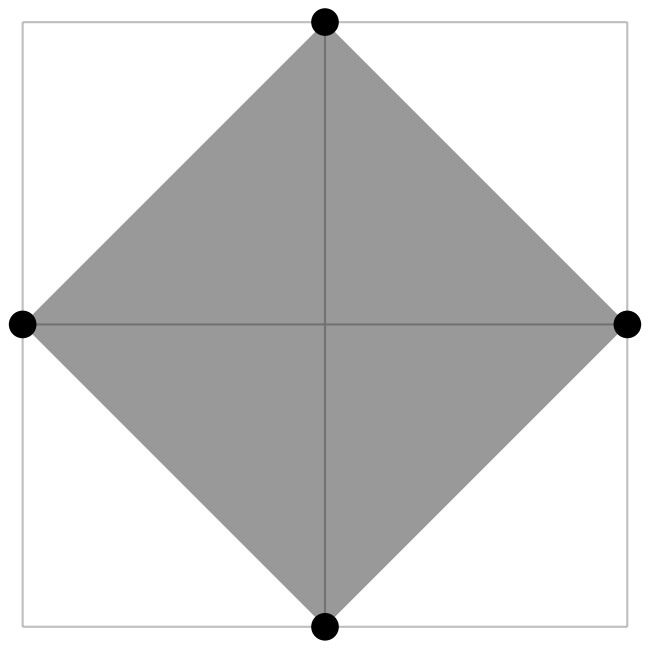

In fact, for each $n\in\mathbb{N}$, $\sqrt{n}$ is a constructible number and the angle $\pi/2$ is a constructible angle, see \cite{Mar1},
so that we can construct a square of area $n\in\mathbb{N}$ in finite steps.
Also, since drawing a square is equivalent to finding corner points of the square,
throughout this paper, we focus on the methods to point the four vertices.
As a consequence, we derive an optimal number of steps for constructing a square of area $n\in\mathbb{N}$ in a grid and 
further, we classify a square of area $n\in\mathbb{N}$ in a grid by an arithmetic property of the number $n\in\mathbb{N}$.

Throughout this paper, 
regarding a single geometric construction step, we follow the standard operational 
definition of ruler-and-compass constructions, see \cite[Definition~2.1]{Mar1}:
drawing a straight line passing through two existing points, or drawing a circle 
centered at an existing point passing through another existing point, 
where any intersections of these lines and circles with each other or with the grid lines are determined without requiring extra steps.
In addition, for convenience, we fix a grid point as the origin $O=(0,0)$ and 
then we assume that every grid point is represented by a pair of integers.

\section{\bf Construction examples of a square of area $3$}
\label{sec2}

Now, we give three examples of geometric construction for finding four vertices of a square of area $3$.
In fact, there are lots of methods and their variations for this.
Here, we point out that we would like to introduce efficient ways based on the contents to be explained in the next section.

\begin{enumerate}

\item[Example 1.]
$3=1^2+1^2+1^2$ \\
In this case, we regard $3$ as the sum of three squares of $1$.
We first draw a circle centered at $(0,-1)$ passing through the point $(1,0)$, which yields the point $(\sqrt{2},-1)$.
Next, by drawing a circle centered at the origin through $(\sqrt{2},-1)$, we simultaneously obtain the locations of $A=(\sqrt{3},0)$ and $B=(0,\sqrt{3})$.
Then, for the line segment $AO$, we draw a circle centered at $A$ with radius $AO=\sqrt{3}$
and we perform the same construction at $B$.
Finally, we let $C$ be the intersection point of circles from $A$ and $B$, which is not the origin. 
The points $O,A,B,$ and $C$ then form the vertices of a square of area $3$.

\begin{figure}[htbp]

  \centering
  \tikzset{
    every node/.style={font=\small},
    grid style/.style={step=1cm, gray!30, very thin},
    axis style/.style={->, >=stealth, thick, black},
    pt/.style={fill=black, draw=black, circle, inner sep=1.5pt}
  }

  \begin{minipage}[b]{0.28\textwidth}
    \centering
    \begin{tikzpicture}[scale=0.85]
      \draw[grid style] (-2.2, -2.5) grid (2.5, 2.2);
      \draw[axis style] (-2.3, 0) -- (2.6, 0) node[right] {$x$};
      \draw[axis style] (0, -2.4) -- (0, 2.3) node[above] {$y$};
      \node[below left] at (0,0) {$O$};


      \draw[thick, black] (0, -1) -- (1, 0);

      \draw[thick, black] ({sqrt(2)*cos(-15)}, {-1 + sqrt(2)*sin(-15)}) arc (-15:90:{sqrt(2)});

      \node[pt] at (0, -1) {};
      \node[below left] at (0, -1) {$(0, -1)$};

      \node[pt] at (1, 0) {};
      \node[above left] at (1, 0) {$(1, 0)$};

      \node[pt] at ({sqrt(2)}, -1) {};
      \node[below right] at ({sqrt(2)}, -1) {$(\sqrt{2}, -1)$};
    \end{tikzpicture}
    \par\small (a) Step 1: Find $(\sqrt{2},-1)$

  \end{minipage}%
  \hfill
  \begin{minipage}[b]{0.28\textwidth}
    \centering
    \begin{tikzpicture}[scale=0.85]
      \draw[grid style] (-2.2, -2.5) grid (2.5, 2.2);
      \draw[axis style] (-2.3, 0) -- (2.6, 0) node[right] {$x$};
      \draw[axis style] (0, -2.4) -- (0, 2.3) node[above] {$y$};
      \node[below left] at (0,0) {$O$};

      \draw[thick, black] ({sqrt(3)*cos(-45)}, {sqrt(3)*sin(-45)}) arc (-45:105:{sqrt(3)});

      \draw[thick, black] (0,0) -- ({sqrt(2)}, -1);

      \node[pt] at (0,0) {};

      \node[pt] at ({sqrt(2)}, -1) {};
      \node[right] at ({sqrt(2)}, -1) {$(\sqrt{2}, -1)$};

      \node[pt] at ({sqrt(3)}, 0) {};
      \node[above right] at ({sqrt(3)}, 0) {$A(\sqrt{3}, 0)$};

      \node[pt] at (0, {sqrt(3)}) {};
      \node[above right] at (0, {sqrt(3)}) {$B(0, \sqrt{3})$};
    \end{tikzpicture}
    \par\small (b) Step 2: Obtain $A(\sqrt{3},0)$ and $B(0,\sqrt{3})$
  \end{minipage}%
  \hfill
  \begin{minipage}[b]{0.28\textwidth}
    \centering
    \begin{tikzpicture}[scale=0.85]
      \draw[grid style] (-2.2, -2.5) grid (2.5, 2.2);
      \draw[axis style] (-2.3, 0) -- (2.6, 0) node[right] {$x$};
      \draw[axis style] (0, -2.4) -- (0, 2.3) node[above] {$y$};
      \node[below left] at (0,0) {$O$};

      \draw[thick, black] ({sqrt(3) + sqrt(3)*cos(75)}, {sqrt(3)*sin(75)}) arc (75:195:{sqrt(3)});

      \draw[thick, black] ({sqrt(3)*cos(-105)}, {sqrt(3) + sqrt(3)*sin(-105)}) arc (-105:15:{sqrt(3)});

      \node[pt] at (0,0) {};

      \node[pt] at ({sqrt(3)}, 0) {};
      \node[below right] at ({sqrt(3)}, 0) {$A(\sqrt{3}, 0)$};

      \node[pt] at (0, {sqrt(3)}) {};
      \node[above left] at (0, {sqrt(3)}) {$B(0, \sqrt{3})$};

      \node[pt] at ({sqrt(3)}, {sqrt(3)}) {};
      \node[above right] at ({sqrt(3)}, {sqrt(3)}) {$C(\sqrt{3}, \sqrt{3})$};
    \end{tikzpicture}
    \par\small (c) Step 3: Construct $C(\sqrt{3},\sqrt{3})$
  \end{minipage}

  \caption{Description of Example 1}
\label{fig2}
\end{figure}
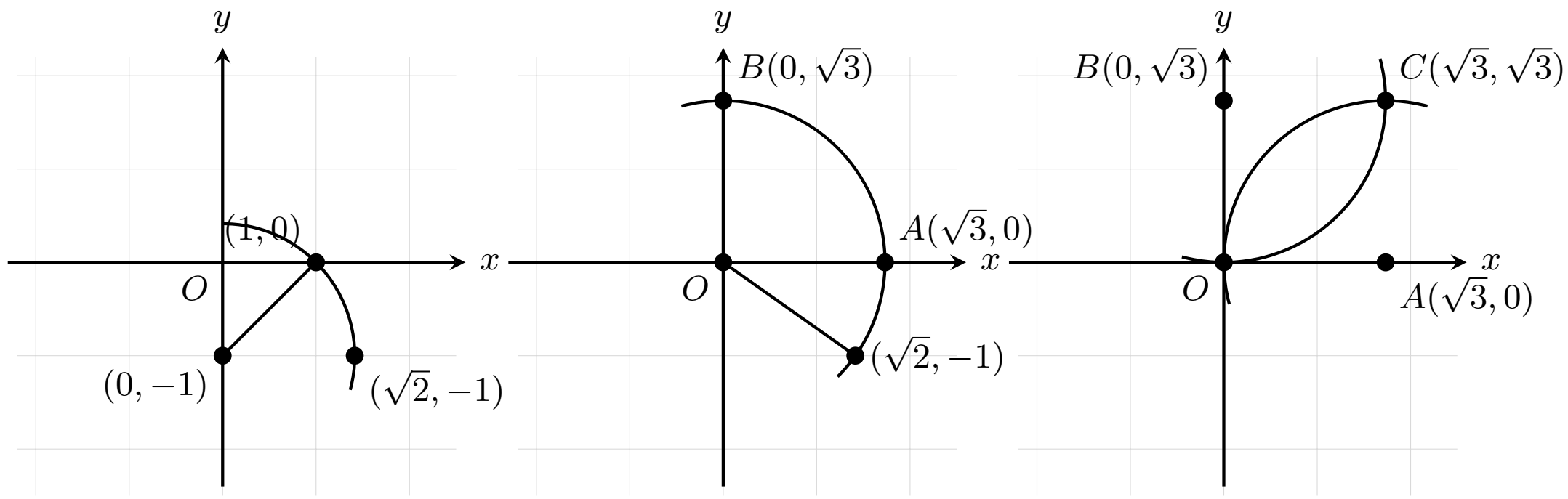

\item[Example 2.]
$3=2^2-1^2$ \\
The Pythagorean theorem $(\sqrt{3})^2+1^2=2^2$ is considered for this case.
We draw a circle centered at $(0,-1)$ with radius $2$, and then
we can find $A=(\sqrt{3},0)$ which is an intersection point of the circle and the line $y=0$.
Next, by drawing a circle centered at the origin passing through $A$, 
we obtain $B=(0,\sqrt{3})$ on the line $x=0$, see Figure \ref{fig3}.
Following the argument of the previous case, then we can determine the final vertex of the square.

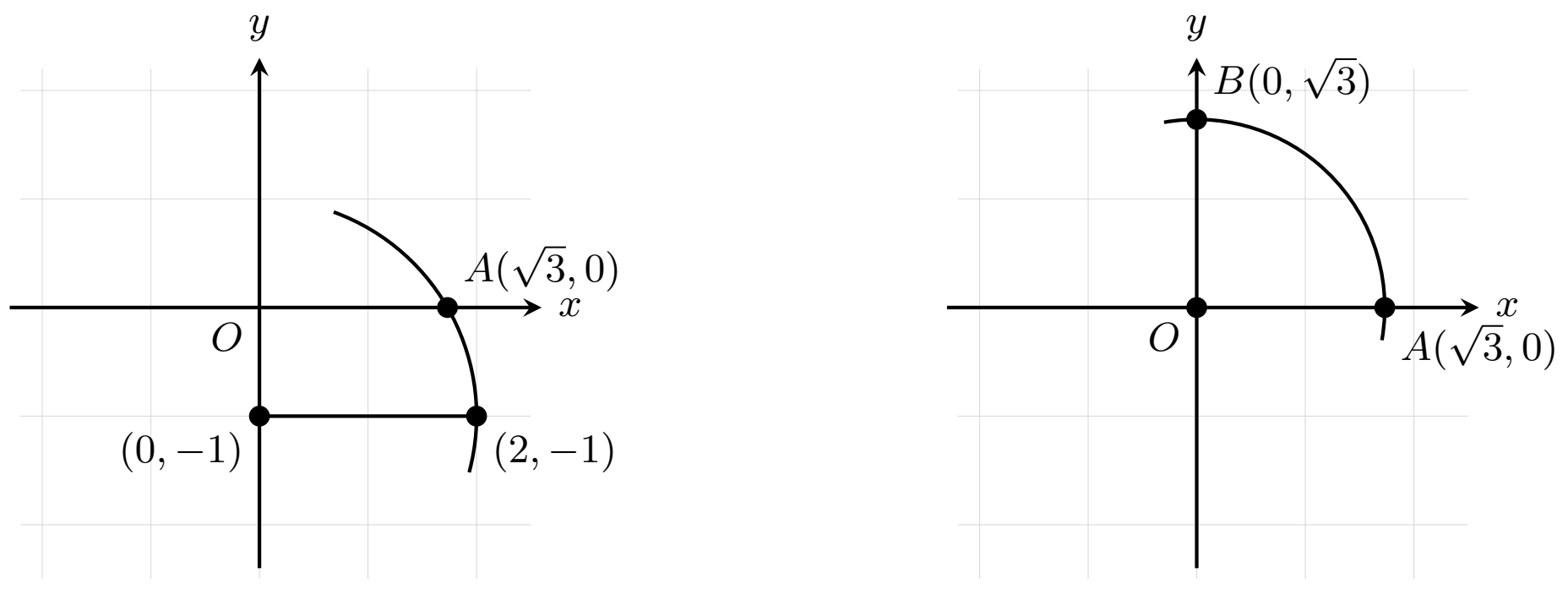
\begin{figure}[htbp]

  \centering
  \tikzset{
    every node/.style={font=\small},
    grid style/.style={step=1cm, gray!30, very thin},
    axis style/.style={->, >=stealth, thick, black},
    pt/.style={fill=black, draw=black, circle, inner sep=1.5pt}
  }

  \begin{minipage}[b]{0.42\textwidth}
    \centering
    \begin{tikzpicture}[scale=0.85]
      \draw[grid style] (-2.2, -2.5) grid (2.5, 2.2);
      \draw[axis style] (-2.3, 0) -- (2.6, 0) node[right] {$x$};
      \draw[axis style] (0, -2.4) -- (0, 2.3) node[above] {$y$};
      \node[below left] at (0,0) {$O$};


      \draw[thick, black] (0, -1) -- (2, -1);

      \draw[thick, black] ({2*cos(-15)}, {-1 + 2*sin(-15)}) arc (-15:70:2);


      \node[pt] at (0, -1) {};
      \node[below left] at (0, -1) {$(0, -1)$};

      \node[pt] at ({sqrt(3)}, 0) {};
      \node[above right] at ({sqrt(3)}, 0) {$A(\sqrt{3}, 0)$};

      \node[pt] at (2, -1) {};
      \node[below right] at (2, -1) {$(2, -1)$};
    \end{tikzpicture}
    \par\small (a) Step 1: Construct $A(\sqrt{3}, 0)$
  \end{minipage}%
  \hfill
  \begin{minipage}[b]{0.42\textwidth}
    \centering
    \begin{tikzpicture}[scale=0.85]
      \draw[grid style] (-2.2, -2.5) grid (2.5, 2.2);
      \draw[axis style] (-2.3, 0) -- (2.6, 0) node[right] {$x$};
      \draw[axis style] (0, -2.4) -- (0, 2.3) node[above] {$y$};
      \node[below left] at (0,0) {$O$};

      \draw[thick, black] ({sqrt(3)*cos(-10)}, {sqrt(3)*sin(-10)}) arc (-10:100:{sqrt(3)});

      \node[pt] at (0,0) {};

      \node[pt] at ({sqrt(3)}, 0) {};
      \node[below right] at ({sqrt(3)}, 0) {$A(\sqrt{3}, 0)$};

      \node[pt] at (0, {sqrt(3)}) {};
      \node[above right] at (0, {sqrt(3)}) {$B(0, \sqrt{3})$};
    \end{tikzpicture}
    \par\small (b) Step 2: Obtain $B(0, \sqrt{3})$
  \end{minipage}

  \caption{First two steps of Example 2}
\label{fig3}
\end{figure}

\item[Example 3.] Using the similarity of right triangles \\
In this method, we use the similarity of right triangles, see Figure \ref{fig4}.
To do this, we place Figure \ref{fig4} in the coordinate system, as shown in Figure \ref{fig5}.
Taking $a=3$, we find $R=(0,\sqrt{3})$, which coincides with $B$ as in Example 1 and Example 2.
Letting $B=R$, we draw a circle centered at the origin passing through $B$, which gives $A=(\sqrt{3},0)$.
Finally, following the final step of Example 1, we determine $C=(\sqrt{3},\sqrt{3})$.

\begin{figure}[htbp]
  \centering
  \tikzset{
    every node/.style={font=\small},
    pt/.style={fill=black, draw=black, circle, inner sep=1.5pt}
  }

  \begin{minipage}[b]{0.5\textwidth}
    \centering
    \begin{tikzpicture}[scale=1.1]
      \coordinate (O) at (0,0);
      \coordinate (P) at (-1,0);
      \coordinate (Q) at (4,0);
      \coordinate (R) at (0,2);

      \draw[thick, black] (P) -- (Q) -- (R) -- cycle;
      \draw[thick, black] (R) -- (O);

      \draw[thin] (0.2, 0) -- (0.2, 0.2) -- (0, 0.2);
      
      \draw[thin] 
        ({-0.2/sqrt(5)}, {2 - 0.4/sqrt(5)}) -- 
        ({-0.2/sqrt(5) + 0.4/sqrt(5)}, {2 - 0.4/sqrt(5) - 0.2/sqrt(5)}) -- 
        ({0.4/sqrt(5)}, {2 - 0.2/sqrt(5)});

      \node[pt] at (P) {};
      \node[pt] at (Q) {};
      \node[pt] at (R) {};
      \node[pt] at (O) {};

      \node[below left] at (P) {$P$};
      \node[below right] at (Q) {$Q$};
      \node[above] at (R) {$R$};
      \node[below] at (O) {$O$};

      \node[below, font=\footnotesize] at (-0.5, 0) {$1$};
      \node[below, font=\footnotesize] at (2, 0) {$a $};
      \node[left, font=\footnotesize] at (0, 1) {$\sqrt{a}$};

    \end{tikzpicture}
    \caption{\mbox{Similarity of right triangles}}
    \label{fig4}
  \end{minipage}
\end{figure}

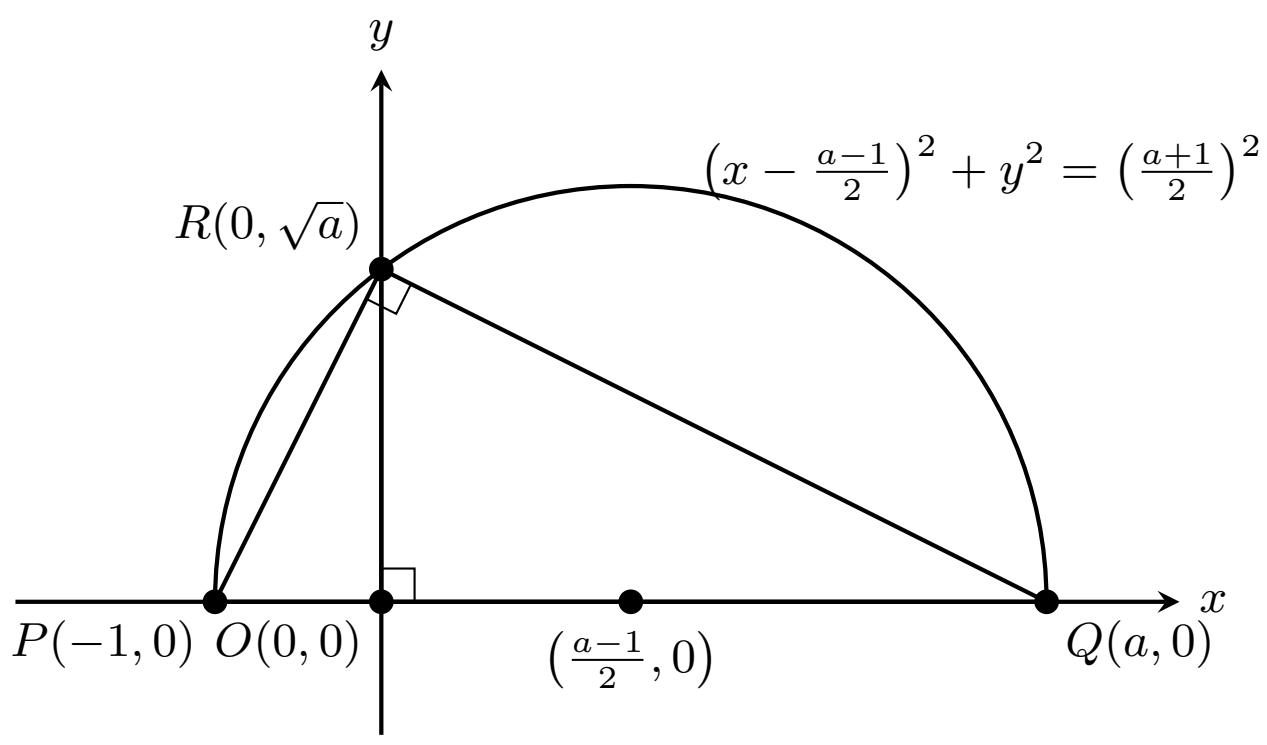
\begin{figure}[htbp]
  \centering
  \tikzset{
    every node/.style={font=\small},
    axis style/.style={->, >=stealth, thick, black},
    pt/.style={fill=black, draw=black, circle, inner sep=1.5pt}
  }

  \begin{minipage}[b]{0.85\textwidth} 
    \centering
    \begin{tikzpicture}[scale=1.1]
      \coordinate (O) at (0,0);
      \coordinate (P) at (-1,0);
      \coordinate (Q) at (4,0);
      \coordinate (R) at (0,2);
      \coordinate (M) at (1.5,0);

      \draw[axis style] (-2.2, 0) -- (4.8, 0) node[right] {$x$};
      \draw[axis style] (0, -0.8) -- (0, 3.2) node[above] {$y$};

      \draw[thick, black] (4,0) arc (0:180:2.5);

      \draw[thick, black] (P) -- (Q) -- (R) -- cycle;
      \draw[thick, black] (R) -- (O);

      \draw[thin] (0.2, 0) -- (0.2, 0.2) -- (0, 0.2);
      
      \draw[thin] 
        ({-0.2/sqrt(5)}, {2 - 0.4/sqrt(5)}) -- 
        ({-0.2/sqrt(5) + 0.4/sqrt(5)}, {2 - 0.4/sqrt(5) - 0.2/sqrt(5)}) -- 
        ({0.4/sqrt(5)}, {2 - 0.2/sqrt(5)});

      \node[pt] at (P) {};
      \node[pt] at (Q) {};
      \node[pt] at (R) {};
      \node[pt] at (O) {};
      \node[pt] at (M) {};

      \node[below left] at (P) {$P(-1,0)$};
      \node[below right] at (Q) {$Q(a,0)$};
      \node[above left] at (R) {$R(0,\sqrt{a})$};
     \node[below left] at (O) {$O(0,0)$};
      \node[below, yshift=-2pt] at (M) {$\left(\frac{a-1}{2}, 0\right)$};

      \node[above right, font=\small] at (1.8, 2.3) 
        {$\left(x - \frac{a-1}{2}\right)^2 + y^2 = \left(\frac{a+1}{2}\right)^2$};

    \end{tikzpicture}
    \caption{Description of Figure \mbox{\ref{fig4}} in the coordinate system}
    \label{fig5}
  \end{minipage}
\end{figure}

\end{enumerate}

From these three examples, we see that the four vertices of a square of area $3$ can be found within four construction steps.
More precisely, at least one step is required to construct the length $\sqrt{3}$.
Then, at least one step is necessary to obtain two line segments of equal length forming a right angle at the origin.
Furthermore, since $C=(\sqrt{3},\sqrt{3})$ does not lie on any grid line, 
at least two additional construction steps are needed to locate $C$.

\section{\bf Some arithmetic properties}
\label{sec3}

In this section, we introduce some arithmetic results for the representation of a positive integer by squares.

\begin{prop}\cite[Theorem 369]{HW1} (Lagrange's four-square theorem)
\label{prop301}
Every positive integer is the sum of four squares.
\end{prop}

Proposition \ref{prop301} means that for any given $n\in\mathbb{N}$, there exist $x,y,z,w\in\mathbb{N}\cup\{ 0 \}$ such that
$$
n=x^2+y^2+z^2+w^2.
$$
Furthermore, the next gives a refinement to the sum of three squares of Proposition \ref{prop301}.

\begin{prop}\cite[Chapter XX]{HW1} (Legendre's three-square theorem)
\label{prop302}
A positive integer $n$ is representable by three squares if and only if 
$n\neq 4^a(8m+7)$ for $a, m\in\mathbb{N}\cup\{ 0 \}$.
\end{prop}

From Proposition \ref{prop302} we see that if $n\in\mathbb{N}$ is not of the form $n=4^a(8m+7)$ for $a, m\in\mathbb{N}\cup\{ 0 \}$,
then there exist $x,y,z\in\mathbb{N}\cup\{ 0 \}$ such that
$$
n=x^2+y^2+z^2.
$$
In order to apply the previous results for the sum of squares to our problem, we introduce the following lemma.

\begin{lem}
\label{lem303}
If a positive integer $n$ is either odd or a multiple of $4$,
then there exist $x,y\in\mathbb{N}\cup\{ 0 \}$ such that
$$
n=x^2-y^2.
$$
\end{lem}

\begin{proof}
If $n=2k+1$ for some nonnegative integer $k$, then
$$
(k+1)^2-k^2=2k+1=n
$$
so that we can take $x=\frac{n+1}{2}$ and $y=\frac{n-1}{2}$.
For the case when $n=4k$ for some $k\in\mathbb{N}$, since
$$
(k+1)^2-(k-1)^2=4k=n,
$$
$x=\frac{n}{4}+1$ and $y=\frac{n}{4}-1$ are the desired numbers.
\end{proof}

According to Proposition \ref{prop302}, it is essential to use Proposition \ref{prop301}
for $n\in\mathbb{N}$ of the form $n=4^a(8m+7)$ for some $a, m\in\mathbb{N}\cup\{ 0 \}$.
However, such integers are either odd or a multiple of 4.
Therefore, we can replace Proposition \ref{prop301} with Lemma \ref{lem303} in our problem under a square grid.

\section{\bf Optimal number of construction steps}
\label{sec4} 

Before presenting our results, we note some observations for the geometric construction on a grid.
\begin{enumerate}

\item
No additional construction steps are needed for a grid point.

\item
At least one step is required for a point lying on a grid line but not a grid point.

\item 
At least two steps are necessary for a point not lying on a grid line.

\item
When we draw a circle, the circle can intersect with several grid lines.

\end{enumerate}

Our first result is the following. This is based on Example 1 and Example 2 of Section 2,
together with Proposition \ref{prop302} and Lemma \ref{lem303} instead of Proposition \ref{prop301} in Section 3.

\begin{thm}
\label{thm401}
For any $n\in\mathbb{N}$, four vertices of a square of area $n$ can be determined within four construction steps.
\end{thm}

\begin{proof}

We divide positive integers into the following cases.

\begin{enumerate}

\item[Case 1.]
$n=x^2+y^2+z^2$ for some $x\in\mathbb{N}$ and $y,z \in\mathbb{N}\cup\{ 0 \}$ \\
Similarly to Example 1 in Section 2,
we first draw a circle centered at $(0,-z)$ passing through the point $(x,y-z)$.
This gives the point $(\sqrt{x^2+y^2},-z)$.
Next, by drawing a circle centered at the origin through $(\sqrt{x^2+y^2},-z)$, 
we obtain the points $A=(\sqrt{x^2+y^2+z^2},0)=(\sqrt{n},0)$ and $B=(0,\sqrt{x^2+y^2+z^2})=(0,\sqrt{n})$ at the same time.
Finally, repeating the procedure described in Example 1 in Section 2 with radius $AO=BO=\sqrt{n}$ determines the final vertex $(\sqrt{n},\sqrt{n})$.
Therefore, four construction steps are needed for such positive integer $n$.

\item[Case 2.]
$n= 4^a(8m+7)$ for some $a, m\in\mathbb{N}\cup\{ 0 \}$ \\
In this case, by Proposition \ref{prop302}, $n$ must be represented by the sum of four squares, so the argument in Case 1 cannot be applied.
On the other hand, since such integer $n$ is either odd or a multiple of $4$,
we apply Lemma \ref{lem303} to $n$ instead of Proposition \ref{prop301}.
Thus, there exist $x,y\in\mathbb{N}$ such that $n=x^2-y^2$ which means the Pythagorean theorem $(\sqrt{n})^2+y^2=x^2$.
With this, similarly to Example 2 in Section 2,
we can find the point $A=(\sqrt{x^2-y^2},0)=(\sqrt{n},0)$ by drawing a circle centered at $(0,-y)$ with radius $x$.
Then we obtain $B=(0,\sqrt{n})$ from the circle centered at the origin with radius $OA=\sqrt{n}$.
Finally, by performing two additional steps as Case 1, we find $C=(\sqrt{n},\sqrt{n})$.
In conclusion, four construction steps are required to this case as well.

\end{enumerate}
\end{proof}

Here, we would like to give some remarks on the proof of Theorem \ref{thm401}.

\begin{rem}
\label{rem402}

In Case 1 of the proof, if $n$ itself is a perfect square, 
then we can take the four vertices as $(0,0)$, $(\sqrt{n},0)$, $(0,\sqrt{n})$ and $(\sqrt{n},\sqrt{n})$ which are grid points,
so no construction steps are required. 
Furthermore, if $n=x^2+y^2$ for some $x, y\in\mathbb{N}$, that is, $n$ is the sum of two squares, 
then $\sqrt{n}$ is obtained directly from a grid point $(x,y)$.
For this case, we do not need any construction steps 
since we can choose the four vertices of a square of area $n$ as $(0,0)$, $(x,y)$, $(y,-x)$ and $(x+y, -x+y)$, all lying on the grid 
and independent of the argument in Case 1.
Therefore, the procedure described in Case 1 is not needed when $n$ is the sum of at most two squares.

We also make the following observation.
For $n\in\mathbb{N}$, when $n$ is neither a perfect square nor the sum of two squares, or when $n= 4^a(8m+7)$ for some $a, m\in\mathbb{N}\cup\{ 0 \}$,
the construction of $\sqrt{n}$ follows the arguments in Theorem \ref{thm401}. 
According to both cases of the proof, to construct the length $\sqrt{n}$,
we see that only the first construction step is required among the four steps,
since the remaining steps merely involve determining the positions of the vertices of the square.

\end{rem}

Now, we further investigate the optimal construction steps for locating four vertices of a square of area $n\in\mathbb{N}$ in the grid.
We will see that the number of steps depends on whether $n$ is even or odd.
Even though the following arguments are applicable to any $n$ only dependent on its parity,
we exclude the cases when either $n$ is a perfect square or the sum of two squares
since such integers do not need any construction steps.
The idea comes from the fact that the length of the diagonal of a square of area $n$ is $\sqrt{2n}$.
First, we state our result for even integers $n\in\mathbb{N}$.

\begin{thm}
\label{thm403}
For any even integer $n\in\mathbb{N}$, four vertices of a square of area $n$ can be determined within two construction steps.
\end{thm}

\begin{proof}
In this case, we let $n=2k$ for some $k\in\mathbb{N}$.
The half length of the diagonal of a square of area $n$ is
$$
\frac{\sqrt{2n}}{2}=\frac{2\sqrt{k}}{2}=\sqrt{k}.
$$
According to the observations at the beginning of Section 4 and Remark \ref{rem402}, for any positive integer $k$, 
we can find a point $(x,y)$ on some grid line with $\sqrt{x^2+y^2}=\sqrt{k}$.
Then, drawing a circle centered at the origin passing through $(x,y)$ yields the four points
$(\sqrt{k},0)$, $(0, \sqrt{k})$, $(-\sqrt{k},0)$, and $(0, -\sqrt{k})$.
These four points form a square whose area is 
$$
\frac{(2\sqrt{k})^2}{2}=\frac{4k}{2}=2k=n.
$$
\end{proof}

\begin{rem}
\label{rem404}
In the proof of Theorem \ref{thm403}, if $k$ is of the form $k=x^2+y^2$ for some $x\in\mathbb{N}$ and $y \in\mathbb{N}\cup\{ 0 \}$ 
then 
$$
(x+y)^2+(x-y)^2=2(x^2+y^2)=2k=n.
$$
This means that such integers are excluded from consideration, so the remaining cases require at least one construction step to obtain the length $\sqrt{k}$.
Moreover, at least one additional step is needed to find the four points forming a square.
For these reasons, the optimal number of steps is two for any even integer $n\in\mathbb{N}$.

\end{rem}

Next, we present our result for odd integers $n\in\mathbb{N}$.
In the following proof, we make use of the similarity of right triangles introduced in Example 3 of Section 2, 
as illustrated in Figure \ref{fig4} and Figure \ref{fig5}.

\begin{thm}
\label{thm405}
For any odd integer $n\in\mathbb{N}$, four vertices of a square of area $n$ can be determined within three construction steps.
\end{thm}

\begin{proof}

Since the case $n=1$ is excluded from consideration, we divide our cases into $n=4k-1$ for some $k\in\mathbb{N}$ and $n=4k+1$ for some $k\in\mathbb{N}$.
Here, we set $a=\frac{n}{2}$ in Figure \ref{fig5}. 

\begin{enumerate}

\item[Case 1.]
$n=4k-1$ for some $k\in\mathbb{N}$ \\
In this case, we first draw the line segment from $(k-1,1)$ to $(k,-3)$, which yields the point $\left(k-1+\frac{1}{4},0\right)=\left(k-\frac{3}{4},0\right)$.
Next, we draw a circle centered at $\left(k-\frac{3}{4},0\right)$ passing through $(-1,0)$. 
Its equation is given by
$$
\left(x-\left(k-\frac{3}{4}\right)\right)^2+y^2=\left(k+\frac{1}{4}\right)^2.
$$
Now, setting $x=0$, we obtain
$$
y^2=\left(k+\frac{1}{4}\right)^2-\left(k-\frac{3}{4}\right)^2=2k-\frac{1}{2}=\frac{4k-1}{2}.
$$
From this, we see that the circle intersects with $x=0$ at the points $\pm\left(0,\sqrt{\frac{4k-1}{2}}\right)=\pm\left(0,\frac{\sqrt{2n}}{2}\right)$.
Since the number $\frac{\sqrt{2n}}{2}$ represents half the length of the diagonal of a square of area $n$, 
by drawing a circle centered at the origin through $\left(0,\frac{\sqrt{2n}}{2}\right)$, 
we obtain the desired four vertices $\left(0,\frac{\sqrt{2n}}{2}\right)$,
$\left(-\frac{\sqrt{2n}}{2},0\right)$, $\left(0,-\frac{\sqrt{2n}}{2}\right)$, and $\left(\frac{\sqrt{2n}}{2},0\right)$.

\item[Case 2.]
$n=4k+1$ for some $k\in\mathbb{N}$ \\
Analogously to Case 1 but with slight modification, for this, we consider the line segment from $(k-1,-3)$ to $(k,1)$
in order to find the point $\left(k-1+\frac{3}{4},0\right)=\left(k-\frac{1}{4},0\right)$.
Using this point as the center, we draw a circle passing through $(-1,0$ which gives the equation
$$
\left(x-\left(k-\frac{1}{4}\right)\right)^2+y^2=\left(k+\frac{3}{4}\right)^2.
$$ 
By taking $x=0$, we have
$$
y^2=\left(k+\frac{3}{4}\right)^2-\left(k-\frac{1}{4}\right)^2=2k+\frac{1}{2}=\frac{4k+1}{2}.
$$
Therefore, this circle intersects with $x=0$ at the points $\pm\left(0,\sqrt{\frac{4k+1}{2}}\right)=\pm\left(0,\frac{\sqrt{2n}}{2}\right)$.
Following the same process as in Case 1, we find the four vertices $\left(0,\frac{\sqrt{2n}}{2}\right)$,
$\left(-\frac{\sqrt{2n}}{2},0\right)$, $\left(0,-\frac{\sqrt{2n}}{2}\right)$, and $\left(\frac{\sqrt{2n}}{2},0\right)$
of the square of area $n$.

\end{enumerate}

\end{proof}

\begin{rem}
\label{rem406}
The construction procedure described in the proof of Theorem \ref{thm405} consists of three steps:
drawing a line segment and constructing two different circles.

To discuss optimality, we consider the case $n=3$.
In this situation, we can easily see that if one vertex of a square of area $3$ lies on a grid point,
then at least one of the remaining vertices does not lie on any grid line.
According to the observations at the beginning of Section 4 and Remark \ref{rem402}, this needs at least three construction steps.
Hence, if we attempt to construct the four vertices within only two steps, none of them can be located at a grid point.

Therefore, we must consider the diagonal of the square.
However, in this case the quantity $\frac{\sqrt{6}}{2}$ cannot be obtained in a single construction step
since constructing the length $\sqrt{6}$ itself requires at least one step and then an additional division process is needed.

From a different point of view, 
all grid points of our grid system are identified with the Euclidean domain of Gaussian integers $\mathbb{Z}[i]=\{a+bi  :  a,b\in\mathbb{Z} \}$, see \cite[Chapter XII]{HW1}.
In $\mathbb{Z}[i]$, to construct the quantity $\frac{\sqrt{6}}{2}$, 
two successive extensions $\mathbb{Z}[i][\frac{1}{2}]$ and $\mathbb{Z}[i][\frac{1}{2},\sqrt{6}]$ are required.
Each algebraic extension corresponds to at least one geometric construction step, 
so that the length $\frac{\sqrt{6}}{2}$ is derived by at least two construction steps.

Consequently, for $n=3$, the least number of construction steps is three and hence
the optimal number of steps is three for any odd integer $n\in\mathbb{N}$.

\end{rem}

Based on the previous arguments, we establish the following classification result.

\begin{thm}
Every positive integer $n\in\mathbb{N}$ can be classified into the following three classes
according to the optimal number of construction steps required to determine the four vertices of a square of area $n$ in a square grid.
\begin{itemize}

\item
$0$-step class:
$n=x^2+y^2$ for some $x\in\mathbb{N}$, $y\in\mathbb{N}\cup\{0\}$.

\item
$2$-step class:
$n$ is an even integer not of the form $n=x^2+y^2$.

\item
$3$-step class:
$n$ is an odd integer not of the form $n=x^2+y^2$.

\end{itemize}
\end{thm}

\begin{proof}

We have already justified the $0$-step class and the $2$-step class in Remark \ref{rem402} and Remark \ref{rem404}, respectively.
For the $3$-step class, 
we first note that if $2n=X^2+Y^2$ with $Z=X+Yi\in\mathbb{Z}[i][\frac{1}{2}]$ for a positive odd integer $n$, then $Z\in\mathbb{Z}[i]$.
Thus, we suppose that there exists a positive odd integer $n$ in this class such that 
$2n=X^2+Y^2$ for some $X\in\mathbb{N}$, $Y\in\mathbb{N}\cup\{0\}$ with $X\geq Y$.
Since $n$ is an odd integer, we can easily see that both $X$ and $Y$ must be odd integers.
This implies that $X-Y$ and $X+Y$ are even integers. Then,
$$
\left(\frac{X-Y}{2}\right)^2+\left(\frac{X+Y}{2}\right)^2=\frac{X^2+Y^2}{2}=\frac{2n}{2}=n
$$
which contradicts the assumption that $n$ is not of the form $n=x^2+y^2$.
This means that for the integer in the $3$-step class, the length $\sqrt{2n}$ cannot be obtained from any grid point.
Therefore, as described in Remark \ref{rem406} for the case $n=3$, 
from $\mathbb{Z}[i]$ two successive extensions $\mathbb{Z}[i][\frac{1}{2}]$ and $\mathbb{Z}[i][\frac{1}{2},\sqrt{2n}]$ are required
for the integer $n$ in this class.

\end{proof}

\section*{\bf Acknowledgment}
The authors would like to thank the students, JunHyeong Cho, Dam Kim, Hyun Lee, Sumin Lee, and Geunmo Park, 
who participated in the Mathematics Program of Science-gifted Education Institute in Kangwon National University.








\bibliographystyle{amsplain}

\end{document}